\input amssym.def
\input amssym
\documentstyle[11pt]{article}
\begin{document}
\newtheorem{thm}{Theorem}[section]
\newtheorem{defin}[thm]{Definition}
\newtheorem{lemma}[thm]{Lemma}
\newtheorem{propo}[thm]{Proposition}
\newtheorem{cor}[thm]{Corollary}
\newtheorem{conj}[thm]{Conjecture}

\centerline{\huge \bf DEFORMATIONS OF}
 
\centerline{\huge \bf(BI)TENSOR CATEGORIES}

\bigskip

\centerline{\parbox{56mm}{Louis Crane and David N. Yetter \\ Department of
Mathematics \\ Kansas State University \\ Manhattan, KS 66506}}
\bigskip
\bigskip

\bigskip

\section{Introduction}

In [1], a new approach was suggested to the construction of four
dimensional Topological Quantum Field Theories (TQFTs), proceeding
from a new algebraic structure called a Hopf category. In [2], it was
argued that a well behaved 4d TQFT would in fact contain such a
category at least formally.

An approach to construction of Hopf categories  was also outlined in
[1], making use of the canonical bases of Lusztig et al [3]. This
proceeded nicely enough, except that there was no natural truncation
of the category corresponding to the case where the deformation
parameter was a root of unity, so all the sums in the tornado 
formula of [1] were
divergent.

This situation is similar to what would have resulted if somebody had
attempted to construct a 3d TQFT before the discovery of quantum
groups. Formal state sums could be written using representations of a
Lie algebra (or its universal enveloping algebra) but they would
diverge. In fact such sums were written in a different context, as
evaluations of spin networks [4]. The discovery of quantum groups made
it possible to obtain finite TQFTs by setting the deformation
parameter equal to a root of unity [5]. 

The key to this progress is the theory of the deformation of Hopf
algebras, as applied to the universal enveloping algebras (UEAs) of simple
Lie algebras. Infinitesimal deformations  
can be classified in terms of a double complex analogous to the
complex which computes the Hochschild cohomology of an algebra.
Certain interesting examples, which lead to global deformations, correspond to Poisson- Lie
algebras, or equivalently to Lie bialgebras, or Manin triples.
Once the interesting infinitesimal deformations of the UEAs were known, it turned
out to be straightforward to extend them to find the quantum groups,
whence the 3d TQFTs followed.

The purpose of this paper is to attempt an analogous procedure for
Hopf categories. We begin by defining a double complex for a ``bialgebra
category,'' whose 3rd cohomology classifies infinitesimal deformations
of the category. Next we apply this complex to the cases of finite
groups and the categorifications of quantum groups produced by
Lusztig [3]. We obtain suggestive preliminary results.

Of course, an infinitesimal deformation is not yet a finite one.
Still less is it a truncation. However, contrary to the folk adage,
lightning tends to strike the same places over and over. The double
complex we construct can also be used to compute the obstructions to
extension of any infinitesimal deformation to a formal series
deformation, so that at least a plausible avenue of research is
opened.

Let us remind the reader of the suggestion that 4d TQFTs may be the
basis for a formulation of the quantum theory of gravity [6]. If this
physical idea is correct, then 4d TQFTs from state sums should exist,
and the program begun in this paper has a good chance of finding them.

In any case, the deformation theory of categories introduced here is
natural, and of intrinsic interest.

The contents of this paper are as follows: chapter 2 describes the
complex which defines the cohomology of a tensor category, and relates
infinitesimal deformations to the third cohomology. Chapter 3 describes the
double complex for a bitensor category, and 
relates it to infinitesimal deformations of
bialgebra categories. Chapter 4 explores the construction of 
infinitesimal deformations in
the most interesting cases.

Let us emphasize that this paper has the purpose of opening a new
direction for research. We pose many more questions than we answer.

\bigskip

\section{Cohomology and Deformations of Tensor Categories}

\bigskip

The deformation theory developed here is very similar in abstract form
to the deformation theory of algebras and bialgebras. Perhaps the not
so categorical reader would do well to study the treatment of that
theory in [5] before reading this chapter. The main formal difference
is that deformations appear in $H^3$ rather than $H^2$. This is
because in a category we deform associators, rather than products, and
similarly for the rest of the structure.

Unfortunately, it will not be practical to make this discussion self
contained. The category theoretic ideas can be found in [7], while the
definition of a bialgebra category is in [1,2,6].

In the following we consider the question of deforming the structure
maps of a tensor category, that is an abelian category $\cal C$ equipped with
a biexact functor $\otimes: {\cal C} \times {\cal C} \rightarrow {\cal C}$
(or equivalently an exact functor $\otimes: {\cal C} \boxtimes {\cal C} 
\rightarrow {\cal C}$, where $\boxtimes$ denotes the universal target
category for biexact functors) which is associative up to a specified
natural isomorphism $\alpha: \otimes (\otimes \times 1) \Longrightarrow
\otimes (1 \times \otimes)$ satisfying the usual Stasheff pentagon.
A tensor category is unital if it is equipped with an object $I$, and
natural isomorphisms $\rho: -\otimes I \Longrightarrow Id_{\cal C}$ and
$\lambda: I\otimes - \Longrightarrow Id_{\cal C}$ satisfying the usual
triangle relation with $\alpha$.

We will consider the case in which the category is $K$-linear for
$K$ some field, usually $\Bbb C$. We denote the category of finite-dimensional
vector-spaces over $K$ by {\bf VECT}.

\begin{defin}
 An {\bf infinitesimal deformation} 
of a $K$-linear tensor category $\cal C$ over
an Artinian local $K$-algebra
$R$ is an $R$-linear tensor category $\tilde{\cal C}$
with the same objects as $\cal C$, but with $Hom_{\tilde{\cal C}} (a,b) =
Hom_{\cal C} (a,b)\otimes_K R$, and composition extended by bilinearity, 
and for which the structure map(s) $\alpha$ ($\rho$ and $\lambda$)
reduce mod $\frak m$ to the structure maps for $\cal C$, where 
$\frak m$ is the maximal ideal of $R$.  A deformation over 
$K[\epsilon]/<\epsilon ^{n+1}>$ is an {\bf $n^{th}$ order deformation}.

Similarly an {\bf $\frak m$-adic deformation} of $\cal C$ over an
$\frak m$-adically complete local $K$-algebra $R$ is an $R$-linear 
tensor category $\tilde{\cal C}$
with the same objects as $\cal C$, but with $Hom_{\tilde{\cal C}} (a,b) =
Hom_{\cal C} (a,b)\widehat{\otimes_K} R$, and composition extended by bilinearity and
continuity, 
and for which the structure map(s) $\alpha$ ($\rho$ and $\lambda$)
reduce mod $\frak m$ to the structure maps for $\cal C$, where 
$\frak m$ is the maximal ideal of $R$. (Here $\widehat{\otimes_K}$ is
the $\frak m$-adic completion of the ordinary tensor product.)  
An $\frak m$-adic deformation over 
$K[[x]]$ is {\bf formal series deformation}.

Two deformations  (in any of the above senses) are {\bf equivalent}
if there exists a monoidal functor, whose underlying functor is
the identity, and whose structure maps reduce mod $\frak m$ to identity
maps from one to the other. 

Finally, if $K = \Bbb C$ (or $\Bbb R$), and all hom-spaces in
$\cal C$ are finite dimensional, a {\bf finite deformation} of $\cal C$ is
a $K$-linear tensor category with the same and maps as $\cal C$, but with
structure maps given by the structure maps of a formal series deformation
evaluated at $x = \xi$ form some $\xi \in K$ such that the formal series
defining all of the structure maps converge at $\xi$.
\end{defin}

Ultimately our interest is in finite deformations, but their study and
construction in general is beyond our present capabilities. In some 
particularly simple cases finite deformations can be constructed directly
(cf. Crane/Yetter [11]).

In the present work, we will confine ourselves to the classification of
first order deformations, and consideration of the obstructions to their
extensions to higher order and formal series deformations.

To accomplish this classification, it is convenient to introduce
a cochain complex 
(over $K$) associated to any $K$-linear tensor category:

First, we fix notation for the totally left and right parenthesized iterates
of $\otimes$ as follows:

\[ \otimes ^n = \otimes(1\boxtimes \otimes)(1\boxtimes 1\boxtimes \otimes)...
(1\boxtimes ... 1 \boxtimes \otimes) \]

\[ ^n \otimes = \otimes(\otimes \boxtimes 1)(\otimes \boxtimes 1\boxtimes 1)...
(\otimes \boxtimes 1 \boxtimes ... 1) \]

\noindent letting $\otimes ^0 = ^0 \otimes = Id_{\cal C}$.

Now, observe that by the $K$-bilinearity of composition,
the collection of natural transformations between any two 
functors targetted at a $K$-linear category forms a $K$-vector space
$Nat[F,G]$. 

We now define the $K$-vector spaces in our complex by

\[ X^n = Nat[^n \otimes, \otimes ^n] \]

Thus elements of $X^n$ have components of the form

\[ f_{A_1,A_2,...,A_n} : (...(A_1\otimes A_2)...\otimes A_n)\longrightarrow
(A_1 \otimes (A_2 \otimes (... \otimes A_n)...)) \]

In order to define the coboundary maps, and in much of what follows, it will
be very convenient to have a notation for a sort of generalized
composition of maps.  To be precise, given some maps $f_1,f_2,...,f_k$
all of whose sources and targets are variously parenthesized tensor
products of the same word of objects, we will denote the composite

\[ a_0 \circ f_1 \circ a_1 \circ f_2 \circ . . . \circ f_k \circ a_k \]

\noindent by $\lceil f_1 f_2 . . . f_k \rceil$  
where $a_0$ is the generalized associator
from the fully left-parenthesized tensor product to the source of $f_1$,
$a_i$ (for $i = 1,...,k-1$) is the generalized associator from the
target of $f_i$ to the source of $f_{i+1}$, and $a_k$ is the generalized
associator from the target of $f_k$ to the fully right-parenthesized 
tensor product.

Except the fact that we need to include $\lceil \;\; \rceil$ 
to obtain well-defined formulas
familiar formulas define the coboundary maps for our complexes:

If $\phi \in X^n$, then $\delta(\phi) \in X^{n+1}$ is defined by

\[ \delta(\phi)_{A_0,...,A_n} = \lceil A_0\otimes \phi_{A_1,...,A_n} \rceil 
+ [\sum_{i=1}^{n-1} (-1)^i \lceil \phi_{A_1,...,A_i\otimes A_{i+1},...,A_n}
\rceil ] + \lceil \phi_{A_0,...,A_{n-1}}\otimes A_n \rceil \]

It follows from the coherence theorem of Mac Lane and the same argument
which show the coboundary in the bar resolution satisfies $\delta ^2 = 0$ that
these coboundaries satisfy $\delta^2 = 0$.  Thus we have a cochain
complex
associated to any $K$-linear tensor category. We denote the cohomology
groups of the complex by ${\frak H}^\bullet({\cal C})$, where the
tensor structure on $\cal C$ is understood. \footnote{In cases where more than
one $K$-linear monoidal structure is being considered on the same category,
it would be necessary to use the notation ${\frak H}^\bullet({\cal C},\otimes,
\alpha)$ to distinguish the structures, since the group depends on the
monoidal structure.}

The significance of this may be found by considering the pentagon relation
for infinitesimal deformations of the associator.

If in the Stasheff pentagon giving the coherence condition on $\alpha$, 
we replace all occurences of $\alpha$ with
occurences of $\alpha + \epsilon a^{(1)}$ (where $\epsilon ^2 = 0$), we 
find that the condition
that the new Stasheff pentagon to commute reduces to

\[ \delta(a^{(1)}) = 0 \]

\noindent where $a^{(1)}$ is considered as an element of $X^3$.

Thus, first order deformations correspond to 3-cocycles in 
our complex.

More, however, is true:  consider now equivalences of first order deformations.
The main (only in the non-unital case) structure map has components of the form

\[ 1_{A\otimes B} + \phi_{A,B} \epsilon : A\otimes B \longrightarrow
A\otimes B \]

\noindent where $\phi$ is some natural endomorphism of $\otimes $.

Now, suppose such a natural transformation defines a monoidal functor
from a first order deformation with associator $\alpha + a \epsilon$
to another with associator $\alpha + b \epsilon$.  Writing out the
hexagon coherence condition for a monoidal functor, and evaluating the
legs then gives

\begin{eqnarray*}
\lefteqn{ \alpha_{A,B,C} +  
\{  a_{A,B,C} + \phi_{A,B\otimes C}(\alpha_{A,B,C}) + 
	(1_A\otimes \phi_{B,C})(\alpha_{A,B,C}) \}\epsilon = } \\
 & & \alpha_{A,B,C} + \{ \alpha_{A,B,C}(\phi_{A\otimes B, C}) +
\alpha_{A,B,C}(\phi_{A,B}\otimes 1_C) + b_{A,B,C} \}\epsilon 
\end{eqnarray*}

Cancelling equal terms, solving for $b_{A,B,C}$, and observing
that the compositions with $\alpha$ are describing the operation of
$\lceil \;\; \rceil$, we find that this
is precisely the condition that

\[ b = a + \delta(\phi) \]

Thus, we have shown:

\begin{thm}
First order deformations of a tensor category ${\cal C}, \otimes, \alpha$
are described by 
3-cocycles in the complex $\{X^n, \delta\}$, and they are classified up 
to equivalence by the cohomology group ${\frak H}^3({\cal C})$.
\end{thm}

Let us now examine the obstructions to extending a first
order deformation to a higher order deformation.

Once again we begin with a commutative Stasheff pentagon, this time with 
legs given by components of $\alpha + a^{(1)}\epsilon$. 

Replacing these with corresponding components of $\alpha + a^{(1)}\eta +
a^{(2)}\eta ^2$ (where $\eta ^3 = 0$), and calculating as before gives us
the condition that

\begin{eqnarray*}
\lefteqn{\delta(a^{(2)}) =} \\
& & \lceil a^{(1)}_{A\otimes B, C, D}a^{(1)}_{A,B,C\otimes D}\rceil - 
\lceil a^{(1)}_{A,B\otimes C,D}(1_A\otimes a^{(1)}_{B,C,D})\rceil \\
& & \mbox{}
- \lceil (a^{(1)}_{A,B,C}\otimes 1_D)(1_A\otimes a^{(1)}_{B,C,D})\rceil -
\lceil (a^{(1)}_{A,B,C}\otimes 1_D)a^{(1)}_{A,B\otimes C,D}\rceil 
\end{eqnarray*}

Thus, the cochain on the right can be regarded as an obstruction to 
the extension to a second order deformation.  It is unclear at this 
writing whether (or under what circumstances) this cochain is
closed. 

In as similar way, the condition needed to extend an $n$-th order
deformation

\[\alpha ^{(n)} = \alpha + a^{(1)}\epsilon + . . . + a^{(n-1)}\epsilon^{n-1} 
\;\;\; (\epsilon^n = 0) \]

to an $n+1$-st order deformation 

\[ \alpha ^{(n+1)} = \alpha + a^{(1)}\eta + . . . + a^{(n)}\eta^{n} 
\;\;\; (\eta^{n+1} = 0) \]

is given by 

\begin{eqnarray*}
\lefteqn{ \delta(a^{(n+1)}) =} \\ 
 & & \sum_{\parbox{.7in}{\scriptsize $i+j = n+1 \\ 1 \leq i,j \leq n$}} 
\lceil a^{(i)}_{A\otimes B, C, D}a^{(j)}_{A,B,C\otimes D}\rceil \\
 & & \mbox{}
- \sum_{\parbox{.7in}{\scriptsize $i+j = n+1 \\ 1 \leq i,j \leq n$}} [
\lceil a^{(i)}_{A,B\otimes C,D}(1_A\otimes a^{(j)}_{B,C,D})\rceil +
\lceil (a^{(i)}_{A,B,C}\otimes 1_D)(1_A\otimes a^{(j)}_{B,C,D})\rceil +
\lceil (a^{(i)}_{A,B,C}\otimes 1_D)a^{(j)}_{A,B\otimes C,D}\rceil ] \\
 & & \mbox{}
- \sum_{\parbox{.7in}{\scriptsize $i+j+k = n+1 \\ 1 \leq i,j,k \leq n$}}
\lceil (a^{(i)}_{A,B,C}\otimes 1_D)
a^{(j)}_{A,B\otimes C,D}(1_A\otimes a^{(k)}_{B,C,D})\rceil
\end{eqnarray*}

\section{Cohomology and Deformations of Bitensor Categories}

A bitensor category is a $K$-linear abelian category  $C$ with two
fundamental structures, a (biexact) tensor product 
$\otimes $ (equivalently, an exact functor 
{\bf $ C \boxtimes C \rightarrow C $ }), which is associative up to a
natural isomorphism $\alpha$ which satisfies the usual Stasheff pentagon,
and a tensor coproduct $\Delta $ 
which is an exact
functor {\bf$ C \rightarrow C \boxtimes C$ } which is coassociative up to
a natural isomorphism $\beta$
which satisfies a dual Stasheff pentagon, and moreover satisfies the condition
that $\Delta$ is a monoidal functor, and $\otimes$ is a cotensor functor
(the dual condition), and the structural transformations are inverse to
each other. We denote the ``coherer'', the structural transformation
for $\Delta$ as a monoidal functor by $\kappa$. A bitensor
category is  biunital when it is equipped with a unit functor $1: {\bf VECT}
\rightarrow C$ and a counit functor  $ \epsilon: C\rightarrow {\bf VECT}$
satisfying the usual triangle, dual triangle and conditions that they 
respect the cotensor and tensor structures, up to mutually inverse 
natural transformations. In the biunital case, we denote the counit
tranformations by $r$ and $l$, and the remaining stuctural transformations
by $\delta$ (counit preserves $\otimes$), $\tau$ (coproduct preserves $I$),
and $\eta$ (counit preserves $I$).

A Hopf category is a biunital bitensor category equipped, moreover, with
an operation on objects, $S$, generalizing dual objects in a suitable sense. 

As in the case of a tensor category, it is the structural 
isomorphisms which we deform, subject to the
coherence axioms. The isomorphisms are natural transformations
between combinations of structural functors, so the terms in finite
order (or formal series) deformations
will live in collections of natural transformations
between the functors, which are vector spaces in the case of
$K$-linear categories
categories and exact functors. 

As we had done for $\otimes$, we fix notation for the 
totally left and right parenthesized iterates
of $\Delta$ as follows:

\[ \Delta ^n = (1 \boxtimes ... 1 \boxtimes \Delta) ... (1\boxtimes \Delta)
\Delta \]

\[ ^n \Delta = (\Delta \boxtimes 1 \boxtimes ... 1) ... (\Delta \boxtimes 1)
\Delta \]

In order for us to place our deformation theory in a cohomological setting,
it will be necessary first to examine the coherence theorem for 
bitensor categories (whether biunital or not).

Fortunately, the structure is given in terms of notions for which
coherence theorems are well known (monoidal categories and monoidal 
functors) or their duals. 

To state it properly, however, we require some preliminaries.  First,
we will restrict our attention to the case where all of our categories are
equivalent as categories without additional structure
to a category $A-{\bf mod}$ for $A$ a finite-dimensional $K$-algebra.
In this case ${\cal C}\boxtimes {\cal D}$ is given by $A\otimes_K B-mod$
when ${\cal C}$ (resp. $\cal D$) is equivalent to $A-{\bf mod}$
(resp. $B-{\bf mod}$).  In this setting, the monoidal bicategory
structure given by $\boxtimes$ has pentagons and triangles which commute
exactly (the structural modifications are identities), so the 
1-categorical coherence theorem of Mac Lane applies, and we may
disregard the parenthesization of iterated $\boxtimes$, and the
intervention of associator and unit functors.
Thus we may use the notation ${\cal C}^{\boxtimes n}$ without
fear of ambiguity.

Second, we must note that if $\cal C$ is a bitensor category, so
is ${\cal C}^{\boxtimes n}$.  The structure functors are given
by applying a ``shuffle'' functor before or after the $\boxtimes$-power
of the corresponding structure functor for $\cal C$.

\begin{thm} {\bf (Coherence Theorem for Bitensor Categories)}
Given two expressions for functors $\Phi, \Phi ^\prime$ 
from an $n$-fold $\boxtimes$-power of a bitensor category to
an $m$-fold $\boxtimes$-power of the same category,
 given in terms of $Id, 
\otimes, \Delta, I, \epsilon, \boxtimes,$ and composition of functors,
(where the structural functors may lie in any $\boxtimes$-power of
$\cal C$), and given two expressions for natural isomorphisms between
these functors in terms of the structural transformations for the
categories, identity transformations,
$\boxtimes$, and the 1- and 2-dimensional compositions
of natural transformations, then in any instantiation of these
expressions by the structures from a particular bitensor category,
the natural isomorphisms named by the two expressions are equal.
\end{thm}

\noindent{\bf proof:}  
First note that as in the corresponding from of Mac Lane's coherence
theorem, we must deal with formal expressions for functors and 
natural transformations to avoid ``coincidental'' compositions.

The proof is reasonably standard:  for any expression
for a functor of the given form, we construct a particular ``canonical''
 expression for
a natural
isomorphism to another such expression for a functor, then show that
given two expressions for functors, and a natural transformation named
by a {\em single instance} of a structural natural isomorphism, 
identity transformations, and 1-dimensional composition of natural
transformations, the diagram of natural transformations formed by
this ``prolongation'' of the structure map
and the two ``canonical'' expressions closes and commutes.

(The ``canonical'' has quotation marks, since it is only once the
theorem is established that we will know that the map named by the
composite is, in fact, canonical.  {\em A priori} it is dependent
upon the construction given.)

Note that this suffices, since 

\begin{enumerate}
\item by the middle-four-interchange law,
any expression for a natural isomorphism of the sort described in
the theorem will factor into a 2-dimensional composition of expressions
of this restricted sort, and
\item any composite of such expressions is then seen to be equal to
the composite of the ``canonical'' expression for the source,
 followed by the inverse
of the ``canonical'' expression for the target.
\end{enumerate}

Our ``canonical'' expressions consists of a composite $c_1$ of instances
of the structure maps $\kappa, \delta, \tau,$ and $\eta$  
to move all occurences of $\Delta$ and $\epsilon$ ``inside'' all
occurences of $\otimes$, and remove all applications of $\Delta$ or
$\epsilon$ to $I$; followed by a composite $c_2$ of  instances of
$\beta, r,$ and $l$ to remove all occurences of $\epsilon$ applied to
a cofactor of $\Delta$ and to completely 
right coassociate all iterated $\Delta$'s; followed by a composite $c_3$
of instances of $\alpha, \rho$ and $\lambda$ to remove all instances
of $I$ tensored with other objects, and completely right associate
all iterated $\otimes$'s.  

Note that we have chosen an order to compose the three constituent
composites, but have not specified the order within each composite.
This is possible because $c_2$ and $c_3$ are each independent
of the order by the coherence theorem of Mac Lane, and its dual, and
the functoriality properties of $\boxtimes$ and the 1-dimensional
composition of natural transformations; while for $c_1$,
the order of application is constrained by the nesting of 
the various functors, but within those constraints, the resulting
composite is independent of the order by virtue of the functoriality
properties of $\boxtimes$ and the 1-dimensional composition of 
natural transformations.

In the circumstances of the theorem, we will let $c_i$ (resp. $c_i^\prime$)
$i=1,2,3$ denote the
components of the ``canonical'' map from $\Phi$ (resp. $\Phi ^\prime$).

We now have three cases

\begin{itemize}
\item[Case 1] The natural isomorphism $f$ from $\Phi$ to $\Phi ^\prime$
 is a prolongation of
$\kappa, \delta, \tau,$ or $\eta$.
\item[Case 2] The natural isomorphism $f$ from $\Phi$ to $\Phi ^\prime$
 is a prolongation of $\beta, r,$ or $l$.
\item [Case 3] The natural isomorphism $f$ from $\Phi$ to $\Phi ^\prime$
 is a prolongation of $\alpha, \rho,$ or $\lambda$.
\end{itemize}

In Case 1, it follows from the same argument that shows that $c_1$ is
well-defined that the targets of $c_1$ and $c_1^\prime$ coincide, and
that $c_1^\prime(f) = c_1$.

In Case 2, by using the functoriality properties of $\boxtimes$ and the
1-dimensional composition of natural transformations, and
the dual of the coherence theorem for monoidal functors, we can construct
a natural isomorphism $f^!$ from the target of $c_1$ to the target of
$c_1^\prime$ such that $f^!$ is a composition of prolongations of 
$\beta$'s, $r$'s, and $l$'s, and $c_1^\prime(f) = f^!(c_1)$.  It then
follows from the same argument that shows $c_2$ is well-defined that the
targets of $c_2$ and $c_2^\prime$ coincide, and $c_2^\prime(f^!) = c_2$.

Finally, for Case 3, by using the functoriality properties of $\boxtimes$
and the 1-dimensional composition, and the coherence theorem for monoidal
functors, we can construct a natural isomorphism $f^!$ from the target 
of $c_1$ to the target of
$c_1^\prime$ such that $f^!$ is a composition of prolongations of
$\alpha$'s, $\rho$'s, and $\lambda$'s, and $c_1^\prime(f) = f^!(c_1)$. 
By using the functoriality properties
of $\boxtimes$ and the 1-dimensional composition, and the naturality 
properties of prolongations of $\beta$, $r$, and $l$, we can construct
a natural isomorphism $f^{!!}$ from the target of 
$c_2$ to the target of $c_2^\prime$
such that $f^{!!}$ is a composition of prolongations of
$\alpha$'s, $\rho$'s, and $\lambda$'s, and $c_2^\prime(f^!) = f^{!!}(c_2)$.
It follows by the same argument that shows $c_3$ is well-defined that the
targets of $c_3$ and $c_3^\prime$ coincide, and $c_3^\prime(f^{!!}) = c_3$.
$\Box$
\bigskip

We shall call an instantiation of expressions of the type given in the
previous theorem a pair of {\em commensurable functors}, and the unique 
natural isomorphism obtained by instantiating an expression of the type 
in the theorem the {\em commensuration}.  Given commensurable 
functors $F$ and $G$,
we will denote the commensuration by $\gamma ^{F,G}$.

Now, observe that $\Delta ^n(^n \otimes)$ and $[\otimes^i]^{\boxtimes j} sh
[^j \Delta]^{\boxtimes i} = \otimes_j^i [^j \Delta]^{\boxtimes i} =
[\otimes^i]^{\boxtimes j} \Delta _i$ are commensurable functors.
Now, given a sequence of natural transformations
$f_1, . . . ,f_n$ such that the source of $f_1$ is commensurable with
$\Delta ^n(^n \otimes)$, and the target of $f_i$ is commensurable with
the source of $f_{i+1}$, and the target of $f_n$ is commensurable with
$[\otimes ^i]^{\boxtimes j} sh
[^j \Delta]^{\boxtimes i}$
let

\[ \lceil f_1,...,f_n \rceil: \Delta ^j(^i \otimes)\Rightarrow 
[\otimes ^i]^{\boxtimes j} sh
[^j \Delta]^{\boxtimes i} \] 

\noindent denote the composition of the given natural transformations
 alternated with the
appropriate commensurations.

For any bitensor category, we can now define a double
complex of vector spaces
 
\[ (X^{\bullet, \ast}, 
d^{\bullet, \ast}:X^{\bullet, \ast}\rightarrow X^{\bullet+1, \ast}, 
\delta^{\bullet, \ast}:X^{\bullet, \ast}\rightarrow 
X^{\bullet, \ast+1}). \]

\noindent Where $X^{ij}$ is the space of natural
transformations between the two (commensurable)
functors $\Delta ^j \; ^i \otimes$ and
$[\otimes ^i]^{\boxtimes j} sh [^j \Delta]^{\boxtimes i}$
from the i-fold to the j-fold
tensor power of {\bf C} to itself, where $sh$ is the ``shuffle functor''
from $[C^{\boxtimes j}]^{\boxtimes i}$ to $[C^{\boxtimes i}]^{\boxtimes j}$. 
(Notice that because our category if $k$-linear, these collections
of natural transformations are $k$-vector spaces.)
And 

\[d(s) = \lceil 
\otimes ^n\boxtimes 1(1\boxtimes s)_{sh(\Delta ^{\boxtimes n}(-)}
\rceil + \left[ \sum_{i=1}^m 
(-1)^i \lceil 1^{i-1}\boxtimes \Delta \boxtimes 1^{m-i} 
\rceil \right] + (-1)^{m+1} \lceil 1\boxtimes \otimes ^n 
(s\boxtimes 1)_{sh(\Delta ^{\boxtimes n}(-)} \rceil  \]

\noindent and

\[ \delta(s) = \lceil
\otimes ^{\boxtimes m}(sh(1\boxtimes s)_{^n\Delta \boxtimes 1 (-)}) \rceil
+ \left[ \sum_{i=1}^n 
(-1)^i \lceil s_{1^{i-1} \boxtimes \otimes \boxtimes 1^{n-i}(-)}
\rceil \right] + (-1)^{n+1}
\lceil \otimes ^{\boxtimes m}(sh(s \boxtimes 1)_{1 \boxtimes ^n\Delta (-)})
\rceil \]

\noindent in each case for $s \in X^{n,m}$

In the case of a biunital bitensor category, we can easily extend our 
complex to
include the values of 0 for i and j, interpreting ${\cal C}^{\boxtimes 0}$
as {\bf VECT},
$\otimes ^0$ and $^0 \otimes$ as the functor $I$, and 
$\Delta ^0$ and $^0 \Delta$ as the functor $\epsilon$.

It follows by a diagram chase from the coherence of the bialgebra
category that $d^2 = \delta ^2 = d \delta + \delta d=0$. Thus we have a
bicomplex whose cohomology can be defined in the usual manner.

\begin{defin} The bicomplex described above is the {\bf basic bicomplex}
of the bitensor category. The total complex of the basic bicomplex, 
indexed by $X^n = \oplus_{i+j = n+1} X^{i,j}$ is the {\bf basic complex}
of the bitensor category.
\end{defin}

\begin{defin}
 The larger bicomplex described above is the {\bf extended bicomplex} of a
biunital bitensor category. 
\end{defin}

Now we note that the three structural natural transformations of a
bialgebra category live in the third diagonal of the basic bicomplex.
Specifically, the associator $ \alpha $ for the tensor product lives in $X^{3,1}$
the coassociator $ \beta$ for the coproduct lives in $X^{1,3}$, and the
``coherer'' $ \kappa $ lives in  $X^{2,2}$. The coherer is the isomorphism in a
bialgebra category that corresponds to the axiom for a bialgebra which
states that $ \Delta (ab)= \Delta (a) \Delta (b) $.

By an infinitesimal deformation of a bialgebra category we mean an
infinitesimal deformation of its structural natural transformations
which satisfies the coherence axioms to first order in the
infinitesimal parameter. this makes sense because natural
transformations are combinations of morphisms, and all the spaces of
morphisms for our spaces are vector spaces.

Concretely we let $ \kappa ' = \kappa + k\epsilon, \alpha '=\alpha
+ a\epsilon, \beta ' =\beta +  b\epsilon$, for $\epsilon ^2 = 0$. 
When we write out the coherence
axioms for the new maps, we find
the only new conditions beyond the coherence of the 
triple $\alpha, \kappa, \beta$ are 

\[d(a) = \delta (a) + d(k) = d(k) + \delta (b) = \delta(b) =0. \]

The deformations of our category (as a bitensor category) 
correspond to cocycles of the basic
complex. Similarly, the equivalence classes of deformations under
infinitesimal monoidal equivalence correspond to cohomology classes.

\begin{thm} The equivalence classes of infinitesimal deformations of
a bialgebra category correspond to classes in the third cohomology of
its basic complex.
\end{thm}

\noindent{\bf proof:} Once it is observed that the structural maps
for a bitensor functor are elements of $X^{1,2}$ and $X^{2,1}$, it
is easy to check (by writing out the hexagon coherence conditions
for monoidal and dual monoidal functors) that a bitensor functor
structure for the identity functor given over $K[\epsilon]/<\epsilon ^2>$
is described by a total 2-cochain which cobounds the difference between
the two bitensor structures (as 3-cochains).

(The fact that the third cohomology group appears here, rather than the
second ala Hochschild, is suggestive in relation to the categorical
ladder picture in TQFT. We know that a TQFT can be constructed from a
finite group plus a cocycle of the group. The cocycle of the group
must be chosen to match the dimension of the TQFT. Thus if a 2-cocycle
of a bialgebra gives rise to a 3d theory, it is plausible that a
3-cocycle of a bialgebra category would generate a 4d theory. All this
raises the question whether there is a classifying space of some sort
for a bialgebra category whose cohomology is related to the cohomology
of our bicomplex.)

This theorem is not very useful in itself, since it does not suggest a
way to find interesting examples of cocycles. However, for a biunital bitensor
category, we can embed the basic bicomplex into the extended
bicomplex. Any element of $X^{0,3}$ on which $\delta$ gives 0, or any
element of $X^{3,0}$ on which $d$ gives 0 can be pushed back into
the basic bicomplex to give a candidate for a deformation. This is
analogous to the process which led to the quantum groups:
the classical r matrix lives in an extended bicomplex, and the
vanishing of  the analog of the Steenrod square of its differential is 
precisely the classical Yang-Baxter
equation. See [5]. (Of course, the classical Yang- Baxter equation was not for
any element of the complex associated to the Hopf algebra, but only to
one of a very special form related to the Lie algebra. At the moment
we do not know an analogous ansatz for the categorified situation. 
)

Thus, we now have an a pair of interesting new equations to
investigate for Hopf categories :

\smallskip

$ d(s) =0, \;\;\; s \in X^{3,0} \;\; (D1)$

\smallskip

or 

\smallskip

$ \delta (t)=0, \;\;\; t \in X^{0,3} \;\; (D2)$

\smallskip

In addition, we can ask about the equation which says that the
infinitesimal defomration contructed from a solution to (D1) or (D2)
can be extended to a second order deformation.

In the bialgebra situation, the combination of these two equations led to
the classical Yang-Baxter equations, in the restricted ansatz.

\section{Searching for Deformations in some Interesting Cases}

\bigskip

A naive reader might suppose that the deformation equations
D1 and D2 are rather disappointing, since they lead to a sort of cohomology of
automorphisms of the identity or counit of the category. However, in
the important cases, the unit and counit are not simple objects, so in fact we
are led into interesting ground.

In the case of the quantum double of the group algebra of a finite
group, the identity is a sum of one ordered pair of group elements for 
each group element. If we
categorify in the natural way, so that each ordered pair of group
elements is a simple object in the category, (see [11]) our 
equation D1 reduces to a cocycle on the group. In effect, we have
reproduced Dijkgraf-Witten theory [9] in the language of deformed Hopf
categories, since the group cocycle for Dijkgraf-Witten theory induces
a finite (and thus) infinitesimal deformation of the Hopf category.

The other interesting case to apply our theory to is the
categorification of the quantized UEA's constructed by Lusztig in his
construction of the canonical bases [10]. In order to get a
construction which worked for the entire QUEA, Lusztig was forced to
replace the identity by a family of projectors corresponding to the
weight lattice. (It must be cautioned that Lusztig only worked things
out explicitly in the
case of SL(2)). Thus the deformation equations translate
into the coboundary equation for the complex for the group cohomology
of the root lattice. This means that possible infinitesimal
deformations of the bialgebra category correspond to 3-forms on the
fundamental torus of the corresponding Lie group. We can see that even at the first order of
deformation theory, our procedure seems to produce something only for
certain Lie algebras- those of rank at least 3. Work is under way to
examine the implications of the second order deformation equations in
this situation.

It seems unlikely that a complete deformation can be found order by
order. Such an approach is too difficult even for bialgebras. Let us
simply cite the fact [12] that it is an open question whether the
vanishing of the obstruction to a second order deformation is always
enough to ensure a deformation to all orders for a bialgebra.
Nevertheless, our preliminary results suggest that deformations may
exist for the bialgebra categories associated to certain special Lie
algebras only. Whether this could bear any relationship to the special
choices of groups which appear in string theory and supergravity is
not clear at the moment, but the possibility cannot be ruled out.

In order to clarify the situation, it will be necessary to find some
global method for producing deformations. As of this writing, we have
only begun to investigate the possibilities. Several lines of thought
suggest themselves:

\begin{enumerate}
\item One could search for a categorified analog of Reshetikhin's proof
that every Lie bialgebra produces a quantum group [13].
In order to attempt this, we need to single out the part of the
bialgebra category of Lusztig corresponding to the Lie algebra itself.
This is rather delicate, since categories do not admit negative
elements, but a way may be found.

\item It is possible to examine special 3-forms on the groups $F_4$ and
$E_6$, related to their constructions from the triality of SO(8).
Perhaps the relationships of these 3-forms with the structure of the
Lie algebras will make it possible to extend the corresponding
cohomology classes of the root lattices to complete deformations of
the corresponding bialgebra categories. If so, the special Lie
algebras for which we can produce bialgebra category deformations will
be physically interesting ones.

\item Lusztig constructed his categories as categories of perverse
sheaves over flag varieties. The flag varieties are known to
have q-deformations in the sense of non-commutative geometry [14].
Perhaps a suitable category of D-modules over the quantum flag algebras
can be constructed.
\end{enumerate}

\bigskip

\section{Conclusions}

\bigskip

Simple Lie groups and Lie algebras are very central constructions in
mathematics. They appear in theoretical physics as the expressions of
symmetry, which is a fundamental principle of that field. It has been
a remarkable recent discovery that the universal enveloping algebras
of Lie algebras, and the function algebras on Lie groups, admit
deformations. This discovery came to mathematics by way of physics.

It is a further remarkable fact that the deformations of the universal
enveloping algebras admit categorifications, i.e. are related to very
special tensor categories.

There is no reason not to try to see if this process goes any farther.
The question whether the categorifications of the deformations can
themselves be deformed is a natural one.

The development in algebra we have outlined has had profound
implications for topology, and at least curious ones for quantum field
theory as well. Perhaps it is puzzling the the categories constructed
by Lusztig do NOT seem to fit into the topological picture surrounding
quantum groups. The direction of work begun in this paper has the
potential of widening the topological picture to include  Lusztig's
categories as well.

Finally, it seems that the relationship between topological
applications of algebraic structures and deformation theory can be
direct. One of us [15] has recently discovered a brief proof of a theorem 
generalizing the well-known result of Birman and Lins [16] that
the coefficients of the HOMFLY and Kauffman polynomials are
Vassiliev invariants. The proof makes a direct connection between the
stratification of the moduli space of embedded curves in $R^3$ and the
deformation theory of {\em braided} tensor categories (cf. Yetter [17]). 
It is plausible to suggest that the deformations
we are attempting to construct may play a similar role.
\newpage

\centerline{\large \bf References}

1. L. Crane and Igor B. Frenkel, Four dimensional topological quantum
field theory, Hopf categories, and the canonical bases, {\em JMP} {\bf 35} 
(10) (1994),
5136ff.

2. L. Crane and D.N. Yetter, On algebraic structures implicit in TQFTs,
{\em JKTR} (to appear).

3. G. Lusztig, {\em Introduction to Quantum Groups}, Birkhauser, Boston, 1993.

4. J.P. Moussouris, Quantum models of spacetime based on recoupling
theory, D.Phil. thesis, Oxford University, 1983.

5. V. Chiari and A. Pressley, {\em A Guide to Quantum Groups}, Cambridge
University Press, Cambridge, 1994.

6. L. Crane, Clock and category, is quantum gravity algebraic? {\em JMP} 
{\bf36} (1995), 6180ff.

7. M. Kapranov and V. Voevodsky, 2-categories and Zamolodchikov's
tetrahedral equation, preprint.

8. S. Maclane, {\em Categories for the Working Mathematician}, Springer
Verlag, Berlin, 1971.

9.  R. Dijkgraf and E. Witten, Topological gauge theories and group
cohomology, {\em CMP} {\bf 129} (1990), 393-429.

10. G. Lusztig, Canonical bases in tensor products, {\em Proc Nat Acad Sci}
{\bf 89} (1992), 8177-8179.

11. L. Crane and D.N. Yetter, Examples of categorification, preprint.

12. M. Gerstenhaber and S. Schack, Algebras, bialgebras quantum groups
and algebraic deformations, in {\em Deformation Theory and Quantum Groups
with Applications to Mathematical Physics} (M. Gerstenhaber and J.
Stasheff, eds.), Contemporary Mathematics 134, AMS, Providence RI, 1992. 

13. N. Reshetiknin, Quantization of Lie bialgebras, {\em Duke Math Journal
Int Math Res Notices}, {\bf 7}, 143-151.

See Also: P. Etingof and D. Kazhdan,  Quantization of Lie bialgebras I, 
e-print q-alg 9506005.

14. V. Lunts and A. Rosenberg, Localization for Quantum groups, preprint.

15. D.N. Yetter, Braided categories and Vassiliev invariants, II, unpublished
lecture, Workshop on Discriminant Loci, Utrecht University, August 26, 1996.

16. J. Birman and X.-S. Lin, Vertex models, quantum groups and Vassiliev's
knot invariants, preprint.

See also:  D. Bar-Natan, On the Vassiliev knot invariants, {\em Topology}
{\bf 34} (2) (1995) 423ff.

17. D.N. Yetter, Framed tangles and a theorem of Deligne on braided
deformations of Tannakian categories, in {\em Deformation Theory 
and Quantum Groups
with Applications to Mathematical Physics} (M. Gerstenhaber and J.
Stasheff, eds.), Contemporary Mathematics 134, AMS, Providence RI, 1992. 

\end{document}